\documentclass[a4paper,oneside]{article}
\usepackage[T1]{fontenc}
\usepackage[frenchb]{babel}
\usepackage{amsmath,amsthm,multirow}
\usepackage{amsfonts}
\usepackage{array}
\usepackage[T1]{fontenc}
\usepackage{lmodern}
\usepackage{graphicx}
\usepackage{textcomp}
\usepackage{morefloats}
\usepackage{color}
\usepackage{times}

\newtheorem{lemme}{Lemme}
\newtheorem{theoreme}{Théorème}
\newtheorem{prop}{Proposition}

\title{Sur le spectre des longueurs des groupes de triangles $(r,p,q)$\\\normalsize{(The length spectra of triangle groups)}}
\author{Emmanuel Philippe\footnote{emmanuel.philippe@ac-aix-marseille.fr}\\Laboratoire Emile Picard\\Université Paul Sabatier Toulouse}
\date{}

\begin{document}
\maketitle
\begin{abstract}
Nous décrivons le début du spectre des longueurs des groupes de
triangles associés à un triangle hyperbolique.\\
\end{abstract}
\begin{center}
\textbf{\small{Abstract}}\\
\small We describe in this report the beginning of the length spectra of fuchsian triangular groups.
\end{center}
\vspace{1cm}
Quand $\Gamma$ est un sous groupe discret du groupe des
isométries directes du demi-plan de Poincaré \textbf{H}, nous étudions la surface
hyperbolique à points coniques $S=\mathbf{H}/\Gamma$. Pour les considérations élémentaires, nous renvoyons à \cite{Ratcliffe} et \cite{Beardon}. Rappelons cependant que le groupe des isométries directes de \textbf{H} contient notamment des éléments appelés hyperboliques qui sont en correspondance avec les géodésiques fermées de la surface et qu'à chaque élément $\gamma$ de ce type, on associe sa longueur
\[ l(\gamma)=\mathrm{Inf}~\{d(x,\gamma x)~;~x\in\mathbf{H}\} \] où $d$ désigne la métrique hyperbolique usuelle sur \textbf{H}.
On considère alors le \textbf{spectre des longueurs} de $\Gamma$, ensemble
des $l(\gamma)$ quand $\gamma$ parcourt les classes de conjugaison
d'hyperboliques de $\Gamma$. Les longueurs sont ordonnées dans l'ordre croissant
en tenant compte de leur multiplicité.

Cet ensemble à été fréquemment étudié depuis l'article de M.Kac
(\cite{Kac}): nous renvoyons par exemple aux travaux de P. Buser,
K-L. Semmler et R. Dianu (\cite{Buser2},\cite{Dianu}) pour des
résultats de rigidité ainsi qu'à ceux, dans un cadre plus numérique,
de R. Vogeler (\cite{Vogeler}). La question essentielle est de
connaître l'information géométrique de $S$ contenue dans le spectre
des longueurs.

Dans l'article \cite{Philippe2}, nous avons déterminé le début du
spectre des longueurs pour les surfaces de genre $0$ à trois points
coniques dont l'un est d'ordre $2$. Cherchant à étendre cette étude
au cas où les trois points coniques sont d'ordres quelconques
$r,p,q$, nous montrons le résultat suivant:\\

\textbf{Théorème A:}\textit{ Soit $3\leq r\leq  p\leq q$ des entiers dont la somme des inverses et strictement inférieure à $1$. Le spectre des
longueurs du groupe de triangle
$\Gamma(r,p,q)$ commence de la manière suivante}\\
$\mathrm{Lsp}~\Gamma(3,3,q)=\{l_1=l_1\dots \}$ \textrm{ pour tout $q\geq 4$ }\\
$\mathrm{Lsp}~\Gamma(3,4,4)=\{l_1=l_1\dots \}$\\
$\mathrm{Lsp}~\Gamma(3,p,q)=\{l_1=l_1<l_3\dots \}$ \textrm{ pour tout $p\geq 4,q\geq5$ }\\
$\mathrm{Lsp}~\Gamma(r,p,q)=\{l_1=l_1\leq l_2=l_2\leq l_3\dots \}$ \textrm{ pour tout $r\geq4$ }\\
\textit{avec $l_1=l_2$ si et seulement si $p=q$ et $l_2=l_3$ si et seulement si $r=p$.}\\

Signalons que
$l_1=2~\mathrm{Argch}~[2\cos\frac{\pi}p\cos\frac{\pi}r+\cos\frac{\pi}q]$,
les autres valeurs étant obtenues par permutation des entiers
$r,p,q$.\\

Nous donnons également à la fin du présent article quelques précisions sur le spectre dans certains cas où $r=2$ et $r=3$, précisions nécessaires si l'objectif est d'identifier un groupe de triangle à l'aide de son spectre des longueurs.

Précisément, nous établissons que\\

$\mathrm{Lsp}~\Gamma(3,4,4)=\{l_1=l_1=\dots<l_3\dots \}$\\
$\mathrm{Lsp}~\Gamma(3,3,q)=\{l_1=l_1\dots<l'_2=l'_2\dots \}$ \textrm{ pour $q=5,6$ }\\
avec
\[l'_2=2~\mathrm{Argch}~[2(\cos\frac{\pi}q)^2+\cos\frac{\pi}q-\frac12]\]
puis que\\

$\mathrm{Lsp}~\Gamma(2,p,\infty)=\{ l_1(2)=l_1(2)<l_2(1,2)\dots \}$ \textrm{ pour $6\leq p\leq 10$}\\
$\mathrm{Lsp}~\Gamma(2,10,10)=\{l_1(2)=l_1(2)<l_2(1,2)\dots\}$\\
avec
\begin{align*}
&l_2(1,2)=2~\textrm{Argch}~[\cos\frac{\pi}q(4(\cos\frac{\pi}p)^2-1)]\\
&l_1(2)=2~\textrm{Argch}~[2\cos\frac{\pi}p\cos\frac{\pi}q]\\
\end{align*}

L'approche utilisée ici est inspirée de l'article de C. Bavard
(\cite{Bavard}) et prolonge celle introduite dans l'article
\cite{Philippe2}. Il s'agit d'introduire une action combinatoire de
$\Gamma$ sur un certain graphe afin de faciliter la description des
éléments hyperboliques. Dans la référence \cite{Philippe1}, on
trouvera également plus de détails calculatoires sur le paragraphe
\ref{s_22} du présent article.\\

\section{Préliminaires}
\subsection{Introduction}
On considère le demi-plan de Poincaré $\mathbf{H}=\{
z\in\mathbf{C}~;~\mathrm{Im}~z >0 \}$ que l'on munit de la distance
hyperbolique $d$.

Si $x,y\in\mathbf{H}$, il existe une unique géodésique $(xy)$ de
$\mathbf{H}$ passant par ces deux points. Le segment géodésique
reliant $x$ à $y$ sera noté $[x,y]$ et la
distance hyperbolique entre $x$ et $y$ sera notée $xy$.

Considérons deux géodésiques $A_1$ et $A_2$ dans $\mathbf{H}$ et choisissons arbitrairement un point sur chacune d'elles: disons $x_1\in
A_1$ et $x_2\in A_2$. Notons $\theta_1$ et $\theta_2$ des mesures des angles orientés $\angle(A_1,x_1x_2)$ et $\angle(x_2x_1,A_2)=\theta_2$. On dispose du lemme suivant:

\begin{lemme}
Les géodésiques $A_1$ et $A_2$ sont disjointes si et seulement si
$\delta(A_2;A_1)>1$ où
\[\delta(A_2;A_1)=|\sin\theta_1\sin\theta_2\cosh
x_1x_2-\cos\theta_1\cos\theta_2| \]
Dans ce cas, on a en outre $d(A_1,A_2)=\cosh \delta(A_2;A_1)$.
\end{lemme}

Si $3\leq r\leq p\leq q$ sont trois entiers, on introduit alors
 \[ \delta(k_1,k_2)=\sin(k_1\pi/p)\sin(k_2\pi/r)\cosh a
-\cos(k_1\pi/p)\cos(k_2\pi/r) \]

où $1\leq k_1\leq p$, $1\leq k_2\leq r$ et
\[ \cosh a=\frac{\cos\frac{\pi}q+\cos\frac{\pi}p\cos\frac{\pi}r}{\sin\frac{\pi}p\sin\frac{\pi}r} \]

Cette quantité mesure la distance hyperbolique entre deux
géodésiques $A_1,A_2$ telles que $x_1x_2=a$, $\theta_2=k_1\pi/p$ et
$\theta_1=k_2\pi/r$.

De manière analogue, on utilisera aussi parfois les notations
suivantes:
 \[ \delta'(k_1,k_2)=\sin(k_1\pi/p)\sin(k_2\pi/q)\cosh b
-\cos(k_1\pi/p)\cos(k_2\pi/q) \]

où $1\leq k_1\leq p$, $1\leq k_2\leq q$ et
\[ \cosh b=\frac{\cos\frac{\pi}r+\cos\frac{\pi}p\cos\frac{\pi}q}{\sin\frac{\pi}p\sin\frac{\pi}q} \]
ainsi que
 \[ \delta''(k_1,k_2)=\sin(k_1\pi/q)\sin(k_2\pi/r)\cosh c
-\cos(k_1\pi/q)\cos(k_2\pi/r) \]

où $1\leq k_1\leq q$, $1\leq k_2\leq r$ et
\[ \cosh c=\frac{\cos\frac{\pi}p+\cos\frac{\pi}q\cos\frac{\pi}r}{\sin\frac{\pi}q\sin\frac{\pi}r} \]

\subsection{Groupe de triangle}
Fixons trois entiers
$r\leq p\leq q$ avec $1/r+1/p+1/q<1$ et considérons un triangle
hyperbolique $T$
d'angles $\pi/r,\pi/p,\pi/q$. On utilisera ici de manière récurrente les notations
\[ X=\cos\frac{\pi}r~;~Y=\cos\frac{\pi}p~;~Z=\cos\frac{\pi}q\]

Signalons que le cercle inscrit à $T$ possède des points de
contact avec les côtés que nous noterons $r^*,p^*,q^*$ avec une
relation évidente ($r^*$ est sur le côté $[p,q]$). Remarquons que
$rp^*=rq^*$, $pq^*=pr^*$ et $qr^*=qp^*$: on notera respectivement
$d_r,d_p$ et $d_q$ ces trois valeurs. Ces quantités s'expriment facilement en fonction de $X,Y,Z$: si on introduit la quantité
 \[ \Delta=X^2+Y^2+Z^2+2XYZ-1 \]
on dispose par exemple des expressions
\[(\sinh d_r)^2=\frac{\Delta}{2(1-X)(1+Y)(1+Z)}~~;~~ \cosh
p^*q^*=\frac{\Delta}{2(1+Y)(1+Z)}+1
\]
D'autre part,le groupe d'isométrie engendré par les
réflexions par rapport aux côtés de $T$ est noté $\Gamma_0(r,p,q)$
et l'ensemble de ses éléments préservant l'orientation est appelé
$\Gamma(r,p,q)$. C'est ce groupe que l'on nommera dans cet article le \textbf{groupe de triangle} associé à $T$.

Le groupe $\Gamma_0(r,p,q)$ admet $T$ comme
domaine fondamental ce qui fournit un pavage $\mathcal{P}_0$ du
demi-plan $\mathbf{H}$ possédant des sommets de valence $r,p$ ou
$q$.
Si l'on retient uniquement les sommets de valence $r,q$ et les
arêtes les reliant, on obtient un nouveau pavage $\mathcal{P}$
constitué de $2p$-gones de côté $c$ avec

\[\cosh c =\displaystyle{\frac{\cos\frac{\pi}p+\cos\frac{\pi}q\cos\frac{\pi}r}{\sin\frac{\pi}q\sin\frac{\pi}r}}
\]

Dans un tel polygone, les
angles aux sommets de valence $q$ valent $2\pi/q$ et ceux aux
sommets de valence $r$ valent $2\pi/r$.\\

Nous expliquons dans les deux prochains paragraphes comment décrire le début du spectre des longueurs du
groupe $\Gamma(r,p,q)$ si $r\geq 3$.\\

\subsection{Action combinatoire du groupe de triangle
}\label{s2_3}

Dans chaque triangle $T$ du pavage $\mathcal{P}_0$, on examine les
cercles hyperboliques centrés aux trois sommets de valences $r,p,q$
et de rayons $d_r,d_p,d_q$: ce sont des cercles tangents aux points
de contact du cercle inscrit à $T$ avec les côtés de $T$. Nous
travaillerons désormais sur l'ensemble de tous les
cercles obtenus de la manière précédente et sur lesquels se situent
les
points du type $r^*,p^*,q^*$.\\

On note dans la suite \textbf{$\mathcal{E}^*$ l'ensemble des points de
type $q^*$}.\\

Sur un cercle dont le centre est un sommet de valence $r$ de
$\mathcal{P}_0$ et dont le rayon est $d_r$, on distingue $r$ points
de type $q^*$ et ceux-ci sont séparés deux à deux par un secteur
angulaire d'angle $2\pi/r$: on dira que deux tels points $x,y$ de
$\mathcal{E}^*$ sont \textbf{voisins} et que l'arc de cercle
$\widetilde{xy}$ les reliant est une \textbf{arête} (de type $r$ si
le cercle est centré en un point de valence $r$ et de type $p$ si le
centre est un point de valence $p$ du pavage initial).

On appelle \textbf{chemin admissible} dans $\mathcal{E}^*$ toute
suite finie d'arêtes consécutives (pour deux éléments distincts de
$\mathcal{E}^*$, il existe au moins un chemin admissible $\beta$
reliant ces deux éléments). La longueur d'un chemin admissible
$\beta$ est le nombre d'arêtes le constituant et sera notée
$L^*(\beta)$.

Si $x$ et $y$ sont deux sommets de $\mathcal{E}^*$, leur
\textbf{distance combinatoire} dans $\mathcal{E}^*$ est notée
$D^*(x,y)$: c'est la longueur minimale d'un chemin $\beta$ reliant
$x$ à $y$.
\[\framebox{$D^*(x,y)=\textrm{Inf}\{ L^*(\beta)~;~\beta\textrm{ reliant $x$ et $y$
 dans $\mathcal{E}^*$}$  \}} \]

Soit $\gamma$ un élément hyperbolique de $\Gamma(r,p,q)$. Celui-ci
laisse stable l'ensemble $\mathcal{E}^*$ et son action sur cet
ensemble est appréciée à l'aide de son \textbf{déplacement
combinatoire minimal}, ou \textbf{niveau}
\[ \framebox{$\lambda^*(\gamma)=\textrm{Inf}~\{ D^*(x,\gamma
x)~;~x\in\mathcal{E}^*$ \}} \]

Décrivons maintenant de manière géométrique les boules combinatoires
dans
$\mathcal{E^*}$.\\

Soit $x_0\in\mathcal{E^*}$ et $n\geq 1$. Si $x\in\mathcal{E}^*$ est
à une distance combinatoire $n$ de $x_0$, il existe un chemin
admissible $\beta$ de longueur $n$ qui relie $x_0$ à $x$. Ce chemin
est formé de $n$ arcs de cercles (ce sont les arêtes) qui sont
chacun de type $r$ ou $p$. Si $a_i\in\{r,p\}$ est le type de la
$i$-ième arête, on dit alors que \textbf{$\beta$ est de type
$(a_1,a_2,\dots,a_n)$}. Si il existe un chemin de type
$(a_1,a_2,\dots,a_n)$ reliant $x_0$ à $x$, le sommet $x$ est
alors dit de type $(a_1,a_2,\dots,a_n)$.\\

Si $n\geq2$ et $x_0\in\mathcal{E^*}$, on définit la quantité
\[ \framebox{$\rho^*(n)=\textrm{Inf}~\{ d(x_0,x)~,~x \textrm{ n'est pas de type $\underbrace{(r,\dots,r)}_{n}$ ni de type $\underbrace{(p,\dots,p)}_{n}$}\}$}\]
Autrement dit, on retire de la sphère combinatoire de rayon $n$ les
quatre points obtenus en tournant $n$ fois consécutivement sur un
des
deux cercles contenant $x_0$ puis on calcule le rayon hyperbolique minimal de la nouvelle sphère combinatoire obtenue.\\

\subsection{Démarche}
Nous détaillons dans ce paragraphe l'approche privilégiée ici
pour décrire toutes les valeurs du spectre des longueurs de
$\Gamma(r,p,q)$ qui sont inférieures à une certaine longueur
choisie $l_0$. Cela repose sur deux lemmes.\\

Le premier lemme (lemme \ref{l_3}) démontre que le pseudo-rayon
hyperbolique $\rho^*$ des boules combinatoires dans $\mathcal{E^*}$
augmente avec le rayon.

Le second lemme (lemme \ref{l_4}) montre qu'un élément hyperbolique
$\gamma$ de $\Gamma(r,p,q)$ ayant une distance de translation petite
possède un déplacement combinatoire minimal $\lambda^*(\gamma)$
borné et que cette borne peut être choisie indépendemment de
$r,p,q$.

On formule alors la propriété \ref{p_1} , qui relie l'action hyperbolique de $\Gamma(r,p,q)$ à l'action combinatoire introduite plus haut.\\

Commençons par établir le
\begin{lemme}\label{l_3}
L'application $\rho^*$ est croissante
\end{lemme}

\textbf{Démonstration:}\\
Constatons que $\rho^*(n)$ n'est définie que pour $n\geq2$. Fixons
$x_0\in\mathcal{E^*}$ et $n\geq3$. Pour tout $x$ n'étant ni de type
$\underbrace{(r,\dots,r)}_{n}$ ni de type
$\underbrace{(p,\dots,p)}_{n}$ vérifiant $D^*(x_0,x)=n$, nous allons
trouver un $y$ tel que

\begin{displaymath}\left\{
\begin{array}{l}
D^{*}(x_0,y)=n-1\\
d(x_0,y)\leq d(x_0,x)\\
y \textrm{ n'est ni de type }\underbrace{(r,\dots,r)}_{n-1} \textrm{
ni de type } \underbrace{(p,\dots,p)}_{n-1}
\end{array}\right.
\end{displaymath}

Pour cela commençons par choisir $\beta$ un chemin admissible minimal reliant
$x_0$ à $x$:
\[\beta=\{\widetilde{x_0x_1},\dots,\widetilde{x_{n-1}x}\} \]

Par minimalité de $\beta$, on a
nécessairement $D^*(x_0,x_{k})=k$ pour tout $1\leq k \leq n-1$.

Si $x_{n-1}$ n'est ni de type
$\underbrace{(r,\dots,r)}_{n-1}$ ni de type
$\underbrace{(p,\dots,p)}_{n-1}$ on choisit donc $y=x_{n-1}$.
Ce point $y$ vérifie bien $yx_0\leq xx_0$: pour s'en convaincre, on peut raisonner par l'absurde en supposons que $\beta$ intersecte au moins deux fois la médiatrice de $[x,y]$ et contredire ainsi la minimalité de $\beta$.

Supposons maintenant que $x_{n-1}$ est de type
$(\underbrace{r,\dots,r}_{n-1})$, le cas où il est de type
$(\underbrace{p,\dots,p}_{n-1})$ se traitant de manière analogue.
On introduit alors $y\in\mathcal{E^*}$ défini par le
chemin $\beta'$ construit à partir de $\beta$ de la manière
suivante:
\begin{displaymath}
\left\{ \begin{array}{l}
\beta'=\{\widetilde{x_0x_1},\dots,\widetilde{x_{n-2}y}\}\\
\angle(x_{n-2}x_{n-3},x_{n-2}y)=\angle(x_{n-1}x_{n-2},x_{n-1}x)
\end{array}\right.
\end{displaymath}
Le sommet $y$ ainsi construit vérifie toutes les conditions souhaitées. $\maltese$\\

Montrons maintenant le deuxième lemme nécessaire:

\begin{lemme}\label{l_4}
Soit $l_0>0$. Soit $\gamma$ un élément hyperbolique de
$\Gamma(r,p,q)$ ayant une distance de translation inférieure ou
égale à $l_0$. Alors il existe au moins un point $x$ de
$\mathcal{E}^*$ tel que $d(x,\gamma x)\leq C^*(l_0)$ avec
\[\framebox{ $\cosh C^*(l_0)=(\cosh c^*)^2(\cosh l_0 -1)+1$} \]
où $c^*=\mathrm{Max}~\{r^*p^*,p^*q^*,q^*r^*\}$.
\end{lemme}

\textbf{Démonstration:}\\
Il s'agit de considérer l'axe de $\gamma$, qui intersecte au moins un
triangle $T^*$ de sommets $r^*,p^*,q^*$ le long de deux côtés. Choisissons alors pour $x$ le sommet de type $q^*$ de ce triangle et on note
$x_0$ son projeté orthogonal sur l'axe. Sa distance à l'axe étant
nécessairement inférieure à $c^*$, il suffit de conclure par un
calcul de trigonométrie hyperbolique. Sachant que $d(x_0,\gamma
x_0)\leq l_0$, on montre alors que
\[ \cosh d(x,\gamma x) \leq (\cosh c^*)^2(\cosh l_0 - 1)+1 \]
ce qui achève la preuve.$\maltese$\\

\begin{prop}\label{p_1}
 Soit $l_0>0$. Soit $\Gamma(r,p,q)$ un groupe de triangle avec $3\leq r\leq p
\leq q$. Il existe un entier $n_0$ tel que
\[\framebox{$l(\gamma)\leq l_0\Rightarrow \lambda^*(\gamma)< n_0$} \]
\end{prop}

\textbf{Démonstration:}\\
Sachant que $\rho^*$ est croissante et non majorée, il existe un
$n_0$ tel que $\rho^*(n_0)>C(l_0)$. Ceci constaté, si $\gamma$ est
un élément hyperbolique de déplacement combinatoire minimal $n\geq
n_0$ alors pour tout élément $x$ de $\mathcal{E}^*$, \[d(x,\gamma
x)\geq \rho^*(n)\geq \rho^*(n_0)> C^*(l_0) \] ce qui achève la
preuve en usant du lemme précédent. $\maltese$\\

Notons que l'on peut trouver $\mathbf{C}(l_0)$ telle que
$C^*(l_0)\leq \mathbf{C}(l_0)$ pour tout $r,p,q$ car $C^*(l_0)$ est
bornée en $r,p,q$. Autrement dit, on peut choisir $n_0$
indépendemment de $r,p,q$.
C'est ce que nous allons faire dans la mise en pratique.\\

\section{Remarques sur le spectre des longueurs}
\subsection{Interprétation géométrique du spectre des
longueurs}\label{s2_4}

Il peut être intéressant de chercher à dessiner les géodésiques
obtenues sur la surface à points coniques $\textbf{H}/
\Gamma(r,p,q)$ à partir des éléments hyperboliques réalisant les
plus petites distances de translation: on cherche
quelle est la forme géométrique des géodésiques les plus courtes sur
les surfaces à points coniques considérées. Rappelons que ces
quotients sont de genre nul et admettent trois points coniques non
équivalents qui correspondent aux trois classes de conjugaison
d'éléments elliptiques présents dans $\Gamma(r,p,q)$.

Si $\gamma=\overline{r}_1\overline{r}_2$ est un élément hyperbolique
de $\Gamma(r,p,q)$ décomposé en produit de deux réflexions d'axes
$A_1,A_2$, sa distance de translation est \[
l(\gamma)=2\textrm{Argch}~d(A_1,A_2) \] et la perpendiculaire
commune aux deux géodésiques $A_1,A_2$ est l'axe de $\gamma$. Cet
axe se projette sur une géodésique fermée du quotient
$\textbf{H}/\Gamma(r,p,q)$ qui a pour longueur $l(\gamma)$. Ainsi,
chaque élément hyperbolique donne naissance à une géodésique fermée
de la surface quotient. Pour obtenir le support géométrique de
cette géodésique, on dessine localement le segment de la
perpendiculaire commune compris entre les géodésiques $A_1,A_2$ et
on double sa longueur.

\subsection{A propos de la multiplicité des longueurs} Soit
$\gamma\in\Gamma(r,p,q)$ un élément hyperbolique. On a clairement
$l(\gamma)=l(\gamma^{-1})$. La question est alors de savoir si cela
augmente la multiplicité de cette longueur dans le spectre.
Cherchant à savoir sous quelles conditions $\gamma$ et $\gamma^{-1}$
sont conjugués dans $\Gamma(r,p,q)$, nous avons établi le fait
suivant:

\begin{lemme}
Soit $\gamma\in\Gamma(r,p,q)$ un élément hyperbolique s'écrivant
$\gamma=\overline{r}_2\overline{r}_1$ avec
$\overline{r}_i\in\Gamma_0(r,p,q)$. Alors $\gamma$ et $\gamma^{-1}$
sont conjugués dans $\Gamma(r,p,q)$ si et seulement si l'axe de
$\gamma$ passe par un point de valence paire du pavage.
\end{lemme}

\textbf{Démonstration:}\\
Supposons donc que $\gamma=\overline{r}_2\overline{r}_1$ et
commençons par établir que $\gamma$ est conjugué à son inverse dans
$\Gamma_0(r,p,q)$ dès que $\overline{r}_2\in\Gamma_0(r,p,q)$, ce qui
est le cas ici. En effet,
\begin{align*}
\overline{r}_2\gamma\overline{r}_2^{-1}&=\overline{r}_2(\overline{r}_2\overline{r}_1)\overline{r}_2^{-1}=\overline{r}_2(\overline{r}_2\overline{r}_1)\overline{r}_2\\
&=\overline{r}_1\overline{r}_2=\gamma^{-1}
\end{align*}

Supposons maintenant que $\gamma^{-1}=x\gamma x^{-1}$ avec
$x\in\Gamma(r,p,q)$. On a alors $x\gamma x^{-1}=\overline{r}_2\gamma
\overline{r_2}^{-1}$ donc $\overline{r}_2 x\gamma x^{-1}=\gamma
\overline{r_2}$ ce qui implique $(\overline{r}_2 x)\gamma
=\gamma(\overline{r_2} x)$. Autrement dit, l'élément
$y=\overline{r}_2x$ commute avec $\gamma$ et laisse donc stable
l'axe de $\gamma$. Comme cet axe est par construction orthogonal à
l'axe de $\overline{r}_2$, on en déduit que $x$ laisse stable l'axe
de $\gamma$. L'élément $x$ étant une isométrie directe, il s'agit soit
d'un hyperbolique de même axe que $\gamma$, soit d'un elliptique d'angle
$\pi$ centré en un point de cet axe.

Supposons que $x$ soit un élément hyperbolique. Comme $x$ et
$\gamma$ se trouvent dans le groupe discret $\Gamma(r,p,q)$, ils
sont tous les deux puissances d'un même élément primitif $\gamma_0$
appartenant à $\Gamma(r,p,q)$. Si $\gamma=\gamma_0^{n}$ et
$x=\gamma_0^{m}$, on aurait alors
\begin{align*}
x\gamma x&= \gamma^{-1}\\
\gamma_0^{m}\gamma_0^{n}\gamma_0^{-m}&=\gamma_0^{-n}\\
\gamma_0^{n}&=\gamma_0^{-n}
\end{align*}
ce qui est impossible si $n\neq0$.

Par conséquent, si $\gamma^{-1}=x\gamma x^{-1}$ avec
$x\in\Gamma(r,p,q)$, alors $x$ est nécessairement un retournement
centré sur l'axe de $\gamma$. La réciproque est immédiate.
$\maltese$\\

Au passage, on constate que si $r,p,$ et $q$ sont impairs, un
élément hyperbolique de $\Gamma(r,p,q)$ n'est jamais conjugué à son
inverse.

\section{Les éléments hyperboliques de niveau $\geq5$}

Pour décrire l'intersection du spectre des longueurs d'un groupe
$\Gamma(r,p,q)$ avec $[0,l_0]$, il suffit, d'après ce qui a été dit
précédemment,
\begin{center}
\fbox{\begin{minipage}{0.8\textwidth}
\begin{itemize}
\item De trouver $n_0$ tel que $\rho^*(n_0)>C^*(l_0)$
\item De décrire les distances de translation de tous les éléments hyperboliques $\gamma$ vérifiant $\lambda^*(\gamma)<n_0$
\item De ranger les distances de translation obtenues dans l'ordre croissant en tenant compte des multiplicités.
\end{itemize}
\end{minipage}}
\end{center}

\subsection{Choix de $l_0$} Soit $T$ un triangle hyperbolique d'angles
$\pi/r,\pi/p,\pi/q$ avec $3\leq r\leq p\leq q$. On choisit
\[
\framebox{$l_0=2~\textrm{Argch}~[2\cos\frac{\pi}p\cos\frac{\pi}q+\cos\frac{\pi}r]$}
\]

Il s'ensuit
\[\cosh C^*(l_0)=\frac{[\Delta+2(1+X)(1+Y)]^2[2(2YZ+X)^2-2]}{4(1+X)^2(1+Y)^2}+1 \]

Ce qui nous intéresse étant le comportement asymptotique de
$C^*(l_0)$, on cherche alors à contrôler cette quantité. On observe
que
\[ \lim_{X,Y,Z\rightarrow 1}{\cosh C^*(l_0)}=37 \]
Reste à trouver $n_0$ tel que $\lim_{X,Y,Z\rightarrow 1}{\cosh
\rho^*(n_0)}>37$ \dots

Nous montrons dans cette section que  $n_0=5$ convient puis nous
décrivons tous les éléments hyperboliques de niveau $1,2,3$ et $4$.

\subsection{Etude de $\rho^*(5)$}\label{s_22} Il s'agit de
démontrer ici que \[ \framebox{$\rho^*(5)> C^*(l_0)$} \] Pour cela,
fixons $x_0\in\mathcal{E}^*$. D'après la définition de $\rho^*$, il
suffit de montrer que tous les points $x~\in\mathcal{E}^*$ tels que
$D^*(x_0,x)=5$ qui ne sont pas du type $(r,r,r,r,r)$ ni du type
$(p,p,p,p,p)$ vérifient $d(x_0,x)> C^*(l_0)$. Pour chaque type fixé,
il suffit de considérer le point $y$ qui minimise la distance
hyperbolique à $x_0$, la détermination de cet élément s'effectuant de manière évidente à l'aide d'arguments géométriques élémentaires.

On applique alors la démarche opératoire suivante:
\begin{center}
\fbox{\begin{minipage}{0.9\textwidth}
\begin{enumerate}
\item On introduit une fonction polynomiale du type
\[ F(X,Y,Z)=4(1+X)^2(1+Y)^2(1+Z)^{\alpha}[\cosh x_0y-\cosh C^*(l_0)] \]
Notre problème se ramenant à prouver que cette fonction est
strictement positive sur un certain domaine. \item On obtient dans un
premier temps un résultat \og asymptotique en $X$ \fg~en majorant
$F$ par une fonction d'une seule variable $X$ qui est strictement
positive si $X$ est assez proche de $1$: on se ramène alors à un
nombre fini de valeurs de $X$ à étudier. \item Pour chaque valeur
$X_0$ de $X$ qui reste à étudier, on dresse le tableau de variation
de $F(X_0,Y,Z)$ en tant que fonction de $Z$ et on montre que, si $Z$
est assez proche de $1$, toutes les fonctions $F(X_0,Y,Z)$ sont
strictement positives: sachant que $Y\leq Z$, il reste alors
éventuellement un nombre fini de valeurs de $Y$ à tester. \item Pour
chaque couple $(X_0,Y_0)$ restant à étudier, on trace la fonction
$F(X_0,Y_0,Z)$ sur $[0,1]$ et on vérifie qu'elle est bien
strictement positive sur le domaine qui nous intéresse.
\end{enumerate}
\end{minipage}}
\end{center}

Les parties calculatoires de cet \og algorithme \fg ont été effectués à l'aide des logiciels \textsf{maple} et \textsf{scilab}.

Remarquons pour conclure qu'il n'a pas été nécessaire d'étudier
tous les types possibles (il y en a 32). Expliquons pourquoi.

Pour calculer $\rho^*(5)$, il faut déterminer, si
$x_0\in\mathcal{E}^*$ est fixé, le minimum de l'ensemble
\[ \Omega=\{ d(x_0,y)~;~\textrm{ $y$ de type $(a_1,\dots,a_5)$ où $a_i\in\{r,p\}$}\}\subset \mathbf{R}^{+} \] Si maintenant $(a_1,a_2,a_3,a_4,a_5)$ est un
quintuplet de $\{r,p\}^5$ et que l'on note
\[ \Omega(a_1,a_2,a_3,a_4,a_5)=\{d(x_0,y)~;~\textrm{ $y$ de type
$(a_1,\dots,a_5)$} \} \] on peut écrire
\[ \Omega=\bigcup_{a_i\in\{r,p\}}{\Omega(a_1,a_2,a_3,a_4,a_5)} \]
Les ensembles définis ci-dessus étant indépendants du
choix de $x_0$ effectué au départ, on en déduit que
\[ \Omega(a_1,a_2,a_3,a_4,a_5)=\Omega(a_5,a_4,a_3,a_2,a_1) \]
Certains cas sont donc redondants et l'on n'est pas tenu d'examiner
les 32 types possibles.

Voici le tableau des correspondances que l'on a utilisées:
\begin{center}
\begin{tabular}{|c|c|}
\hline \textrm{Cas étudiés} & \textrm{Cas associés}\\
\hline
$(r,r,r,r,p)$ \textrm{ et } $(p,p,p,p,r)$ & $(p,r,r,r,r)$ \textrm{ et } $(r,p,p,p,p)$\\
$(r,r,r,p,p)$ \textrm{ et } $(p,p,p,r,r)$ & $(p,p,r,r,r)$ \textrm{ et } $(r,r,p,p,p)$\\
$(r,p,p,p,r)$ \textrm{ et } $(p,r,r,r,p)$ & $\emptyset$ \\
$(r,p,p,r,r)$ \textrm{ et } $(p,r,r,p,p)$ & $(r,r,p,p,r)$ \textrm{ et } $(p,p,r,r,p)$\\
$(r,r,r,p,r)$ \textrm{ et } $(p,p,p,r,p)$ & $(r,p,r,r,r)$ \textrm{ et } $(p,r,p,p,p)$\\
$(r,r,p,r,r)$ \textrm{ et } $(p,p,r,p,p)$ & $\emptyset$\\
$(r,p,r,r,p)$ \textrm{ et } $(p,r,p,p,r)$ & $(p,r,r,p,r)$ \textrm{ et } $(r,p,p,r,p)$\\
$(r,r,p,r,p)$ \textrm{ et } $(p,p,r,p,r)$ & $(p,r,p,r,r)$ \textrm{ et } $(r,p,r,p,p)$\\
$(r,p,r,p,r)$ \textrm{ et } $(p,r,p,r,p)$ & $\emptyset$\\
\hline
\end{tabular}
\end{center}

Finalement, nous pouvons établir suite à cette étude la\\

\begin{prop}\label{p_12}
Soit $\Gamma(r,p,q)$ un groupe de triangle avec $r\geq 3$ et \[
(r,p,q)\neq (3,4,4)~;~(4,4,4)~;~(5,5,5) \] Alors tout élément
hyperbolique $\gamma$ de ce groupe vérifiant $\lambda^*(\gamma)\geq
5$ possède une distance de translation strictement supérieure à
$l_0$.
\end{prop}

\textbf{Démonstration:}\\
Il suffit d'appliquer la propriété \ref{p_1} après avoir constaté
que $\rho^*(5)>C^*(l_0)$, ce qui a fait l'objet du présent paragraphe. $\maltese$\\

Remarquons que les trois cas exceptionnels correspondent tous à une situation où la fonction $F$ associée au type $(r,p,r,p,r)$ se trouve être négative ou nulle, ce qui n'est pas le cas pour les autres types.\\

Pour décrire les valeurs du spectre des longueurs de $\Gamma(r,p,q)$
qui se trouvent dans l'intervalle $[0,l_0]$, il suffit donc
d'étudier maintenant les éléments hyperboliques $\gamma$ tels que
$\lambda^*(\gamma)\leq4$.

\section{Les éléments hyperboliques de niveau $\leq 4$}\label{s4_3}

On suppose dans cette section que $\gamma$ est un élément
hyperbolique du groupe $\Gamma(r,p,q)$ ayant un déplacement combinatoire minimal inférieur ou égal à $4$, c'est à
dire $\lambda^*(\gamma)\leq 4$. Nous voulons calculer les distances
de translation de tous les éléments hyperboliques du groupe ayant
cette propriété: expliquons comment s'y prendre.

Si $\gamma$ possède un déplacement combinatoire minimal $n$, il
existe un sommet $x\in\mathcal{E}^*$ tel que $D^*(x,\gamma x)=n$. On
choisit alors arbitrairement $T$, un des deux triangles admettant
$x$ comme point de type $q^*$. Le sommet $\gamma x$ est lui aussi un
point de type $q^*$ pour deux triangles possibles. Sachant que
$\gamma$ est une isométrie directe, un seul de ces deux triangles
peut convenir pour $\gamma(T)$ et le couple $(T,\gamma(T))$
caractérise $\gamma$ dans tout le groupe des isométries de
$\mathbf{H}$. On peut ainsi déterminer explicitement $\gamma$ en
trouvant un élément du groupe $\Gamma(r,p,q)$ transformant $T$ en
$\gamma(T)$. On calcule alors sa distance de translation s'il s'agit
d'un hyperbolique.\\

Il est vivement conseillé au lecteur de s'appuyer sur une figure locale du pavage pour se persuader des résultats énoncés dans cette section.

\subsection{Les éléments de déplacement
combinatoire $1$} Il existe dans ce cas un $x\in\mathcal{E^*}$ tel
que $D^*(x,\gamma x)=1$. Choisissons un chemin $\beta$ de longueur
minimale reliant $x$ à $y=\gamma x$: $\beta$ est de type $(r)$ ou de
type $(p)$ et la seule isométrie directe transformant $x$ en $\gamma
x$ est un elliptique de centre $r$ (resp. $p$) et d'angle $\pm
2\pi/r$ (resp. $\pm 2\pi/p$). Ainsi, il n'existe pas d'élément
hyperbolique de déplacement combinatoire minimal égal à $1$ dans
$\Gamma(r,p,q)$.

\subsection{Les élements de déplacement
combinatoire $2$}

Soit $x$ tel que $D^*(x,\gamma x)=2$ et $\beta$ un chemin minimal
reliant $x$ à $\gamma x$. Ce chemin peut être de type
$(r,r),(p,p),(r,p)$ ou $(p,r)$. Nous examinons uniquement ici les cas  $(p,p)$ et
$(r,p)$, les deux autres situations s'en déduisant en considérant
$\gamma^{-1}$ ou en permuttant le rôle de $r$ et $p$. Pour chacun des types étudiés, on commence par fixer
l'orientation de la première arête du chemin $\beta$, l'autre
orientation donnant un élément conjugué dans $\Gamma_0(r,p,q)$ (via
la réflexion d'axe $(rp)$ où $x\in (rp)$) à celui obtenu: les
distances de
translation obtenues sont identiques.\\

Si $\beta$ est de type $(r,r)$, les seules isométries directes
transformant $x$ en $\gamma x$ sont des éléments elliptiques
d'angle $\pm 4\pi/r$.\\

Si le chemin $\beta$ est de type $(r,p)$, il existe quatre solutions
pour $\gamma x$, chacune donnant naissance à une isométrie directe.

La première est associée à un
élément elliptique d'angle $-2\pi/q$.

La deuxième donne un élément hyperbolique de distance
de translation
\[l(\gamma)=2~\textrm{Argch}~[\delta(p-1,1)]\] Pour les deux autres
isométries directes annoncées, elles sont obtenues en conjuguant par
une réflexion. On trouve donc un elliptique d'angle
$-2\pi/q$ et un élément hyperbolique de même distance de translation
que le précédent
\[l(\gamma)=2~\textrm{Argch}~[\delta(1,r-1)] \]

Pour un chemin de $\beta$ de type $(p,p)$, on trouve de même deux
éléments elliptiques. Si $\beta$ est de type $(p,r)$, on trouve deux
éléments hyperboliques de même distance de translation que ceux
trouvés ci-dessus et qui correspondent aux transformations inverses
$\gamma^{-1}$.

\subsection{Les élements de déplacement
combinatoire $3$}

Soit $x\in\mathcal{E}^*$ tel que $D^*(x,\gamma x)=3$ et $\beta$ un
chemin minimal reliant $x$ à $\gamma x$. Pour les raisons déjà
invoquées, nous examinons ici en priorité les cas où $\beta$ est
d'un des types suivants: $(r,r,r),(r,p,r)$ et $(r,r,p)$. Pour
obtenir tous les cas, il suffit de permuter le rôle de $r$ et $p$
dans les décompositions obtenues ou bien de constater que si
$\gamma$ est associé à un chemin de type $(a_1,a_2,a_3)$, alors
$\gamma^{-1}$ est associé à un chemin de type $(a_3,a_2,a_1)$. Voici
le tableau des
correspondances utilisées dans ce paragraphe.\\
\begin{center}
\begin{tabular}{|c|c|}
\hline \textrm{Cas étudié pour $\beta$} & \textrm{Cas associés}\\
\hline $(r,r,r)$ & $(p,p,p)$\\
$(r,p,r)$ & $(p,r,p)$ \\
$(r,r,p)$ & $(p,p,r)$,  $(p,r,r)$ \textrm{ et } $(r,p,p)$\\
\hline
\end{tabular}\\
\end{center}
\vspace{0.5cm}

Si $\beta$ est de type $(r,r,r)$, alors les seules isométries
directes transformant $x$ en $\gamma x$ sont des éléments
elliptiques centrés en un point de valence $r$ et d'angle $\pm
6\pi/r$. Pour le cas $(p,p,p)$, on trouve des éléments elliptiques
centrés en un point de valence $p$ et d'angle
$\pm 6\pi/p$.\\

Dans le cas $(r,p,r)$, seize isométries transforment $x$ en $\gamma
x$, dont huit seulement sont directes.

On trouve deux éléments elliptiques et deux éléments hyperboliques de distances de translation
\[2~\textrm{Argch}~[\delta(1,r-2)]~~\textrm{ et }~~2~\textrm{Argch}~[\delta(p-1,r-2)]\]

Parmi les quatre isométries qui correspondent au cas où la première arête
de $\beta$ est orientée dans l'autre sens on retrouve les
mêmes distances de translation, et les deux éléments hyperboliques qui
correspondent sont conjugués dans $\Gamma(r,p,q)$ aux inverses de
ceux étudiés ci-dessus.

Si maintenant $\beta$ est de type $(p,r,p)$ et si
$\delta(2,r-1),\delta(2,1)>1$, on trouve quatre transformations
hyperboliques inverses deux à deux (à conjugaison près) comme dans
le cas $(r,p,r)$. Les distances de translation de ces éléments sont
\[2~\textrm{Argch}~[\delta(2,r-1)]~~\textrm{et}~~2~\textrm{Argch}~[\delta(2,1)]\]

Dans le cas $(r,r,p)$, il y a huit isométries possibles de ce
type dont quatre directes, toutes hyperboliques si $\delta(1,2)>1$
et $\delta(p-1,2)>1$. Les longueurs recensées sont alors
\[2~\textrm{Argch}~[\delta(1,2)]~~\textrm{ et }~~2~\textrm{Argch}~[\delta(p-1,2)]\]
Comme dans le cas précédent, les deux autres éléments hyperboliques,
qui s'obtiennent en renversant l'orientation de la première arête de
$\beta$, sont conjugués à ceux obtenus ici par une réflexion: ils
possèdent donc les mêmes distances de translation.\\

Si $\beta$ est de type $(p,r,r)$, on trouve, à conjugaison près dans
$\Gamma(r,p,q)$, les tranformations inverses de celles trouvées au
cas $(r,r,p)$: les distances de translation sont identiques. Pour
les cas $(p,p,r)$ et $(r,p,p)$, les hyperboliques trouvés dans le
second cas correspondent aux inverses de ceux du premier cas. Ils
ont, quand $\delta(2,1)>1$ et que $\delta(p-2,1)>1$, les distances
de translation suivantes, obtenues en permuttant les rôles de $r$ et
de $p$ dans l'étude précédente:
\[2~\textrm{Argch}~[\delta(2,1)]~~\textrm{et}~~2~\textrm{Argch}~[\delta(p-2,1)]\]

\subsection{Les éléments de déplacement combinatoire $4$}

Soit $x\in\mathcal{E}^*$ tel que $D^*(x,\gamma x)=4$ et $\beta$
minimal reliant $x$ à $\gamma x$. Nous examinons ici en détail les
cas où $\beta$ est d'un des types suivants:
\[(r,r,r,r),(r,r,r,p),(r,r,p,r),(r,p,r,p),(r,r,p,p)\textrm{ ou } (r,p,p,r)\]Les autres cas s'en déduisent comme précédemment en utilisant ce tableau des correspondances:\\
\begin{center}
\begin{tabular}{|c|c|}
\hline \textrm{Cas étudié pour $\beta$} & \textrm{Cas associés}\\
\hline $(r,r,r,r)$ & $(p,p,p,p)$\\
$(r,r,r,p)$ & $(p,p,p,r)$, $(r,p,p,p)$ \textrm{ et } $(p,r,r,r)$  \\
$(r,r,p,r)$ & $(p,p,r,p)$, $(r,p,r,r)$ \textrm{ et } $(p,r,p,p)$\\
$(r,p,r,p)$ & $(p,r,p,r)$\\
$(r,r,p,p)$ & $(p,p,r,r)$\\
$(r,p,p,r)$ & $(p,r,r,p)$\\
\hline
\end{tabular}
\end{center}
\vspace{0.5cm}

Si $\beta$ est de type $(r,r,r,r)$ avec une première arête orientée
dans le sens direct, une seule isométrie directe transforme $x$ en
$\gamma x$. C'est un elliptique d'angle $8\pi/r$. Si la première
arête est orientée dans le sens indirect, l'angle est $-8\pi/r$.

Si $\beta$ est de type $(p,p,p,p)$, on trouve des
elliptiques d'angle $\pm8\pi/p$.\\

Si $\beta$ est de type $(r,r,r,p)$, il y a huit isométries à étudier
dont quatre isométries directes, toutes hyperboliques. Elles se répartissent par couple en fonction des deux longueurs suivantes:
\[2~\textrm{Argch}~[\delta(1,3)]~~\textrm{ et }~~2~\textrm{Argch}~[\delta(p-1,3)]\]

Si le chemin $\beta$ est de type $(p,p,p,r)$, on trouve, en
permuttant les rôles de $r$ et de $p$ dans l'étude précédente,
quatre isométries directes dont les distances de translation sont
(si $\delta(p-3,1)>1$ et $\delta(p-3,r-1)>1$):
\[2~\textrm{Argch}~[\delta(p-3,1)]
~~\textrm{et}~~2~\textrm{Argch}~[\delta(p-3,r-1)]\]

Si $\beta$ est de type $(r,p,p,p)$ ou $(p,r,r,r)$, l'isométrie
$\gamma^{-1}$, qui possède la même distance de translation que
$\gamma$, figure parmi les cas étudiés ci-dessus: on retrouve les
même
valeurs.\\

Dans le cas $(r,r,p,r)$, on trouve seize isométries à étudier pami lesquelles
huit isométries directes. On y décèle six hyperboliques répartis par couple selon les longueurs suivantes

\[2~\textrm{Argch}~[\delta(p-1,r-3)]~~\textrm{;}~~2~\textrm{Argch}~[\delta(1,r-3)]~~\textrm{ et }~~2~\textrm{Argch}~[\delta(1,r-1)]\]

Si $\beta$ est du type $(p,p,r,p)$, on obtient en permutant les
rôles de $r$ et $p$ six éléments hyperboliques ayant les distances
de translation suivantes:
\[2~\textrm{Argch}~[\delta(3,1)]~~\textrm{;}~~l(\gamma)=2~\textrm{Argch}~[\delta(3,r-1)]~~\textrm{ et }~~2~\textrm{Argch}~[\delta(1,r-1)]\]

Si $\beta$ est du type $(r,p,r,r)$ ou $(p,r,p,p)$, la transformation
$\gamma^{-1}$ figure dans les cas précédents et on
retrouve des distances de translation déjà rencontrées.\\

Dans le cas $(r,p,r,p)$, il y a 32 isométries dont 16 isométries
directes. Parmi les isométries directes, outre un élément elliptique, on trouve deux hyperboliques de longueur
\[2~\textrm{Argch}~[\delta'(2,2)]\]
ainsi que quatre hyperboliques, répartis en couples, qui sont obtenus cette fois en composant deux elliptiques de centres distincts. Les distances obtenues sont alors
\[
2~\textrm{Argch}~[4\cos\frac{\pi}p\cos\frac{\pi}r\cos\frac{\pi}q+2(\cos\frac{\pi}r)^2+2(\cos\frac{\pi}q)^2+2(\cos\frac{\pi}p)^2-1]
\]
et
\[
2~\textrm{Argch}~[4\cos\frac{\pi}p\cos\frac{\pi}r\cos\frac{\pi}q+2(\cos\frac{\pi}r)^2+2(\cos\frac{\pi}q)^2-1]
\]
Ces éléments
sont les deux seuls obtenus ici ne s'écrivant pas comme la composée de deux réflexions ayant pour axes des géodésiques du pavage initial.\\
On trouve également dans ce cas un hyperbolique de longueur
\[4~\textrm{Argch}~[\delta(1,r-1)]\]

Les isométries directes restantes, qui s'obtiennent avec une arête
de départ de sens opposé pour $\beta$, possèdent les mêmes distances
de translation car elles sont conjuguées à celles étudiées ici par
une réflexion.

Si $\beta$ est du type $(p,r,p,r)$, on retrouve des distances de
translation déjà rencontrées dans le cas précédent car si $\gamma$ est de associé au type $(p,r,p,r)$, l'élément $\gamma^{-1}$ est associée au type $(r,p,r,p)$.\\

Examinons maintenant le cas $(r,r,p,p)$. Il y a cette fois huit isométries  à
étudier dont quatre directes, toutes hyperboliques. Elles se répartissent par couple (étant conjuguées deux à deux par une réflexion) et donnent naissance aux longueurs suivantes:
\[2~\textrm{Argch}~[\delta(2,2)]~~\textrm{ et }~~2~\textrm{Argch}~[\delta(p-2,2)]\]

Si $\beta$ est du type $(p,p,r,r)$, la transformations $\gamma^{-1}$
est du type précédent donc aucune nouvelle distance
de translation n'apparaît.\\

Finalement, étudions le cas $(r,p,p,r)$: on dispose de seize isométries dont
huit isométries directes. Outre quatre éléments elliptiques, on y trouve quatre éléments hyperboliques, répartis en couples car conjuguées via une réflexion, qui correspondent aux longueurs
\[2~\textrm{Argch}~[\delta(p-2,r-2)]~~\textrm{ et }~~2~\textrm{Argch}~[\delta(p-2,2)]\]

Si $\beta$ est du type $(p,r,r,p)$, l'isométrie $\gamma^{-1}$ se
trouve associée au cas précédent: les distances de translation
sont donc les mêmes.\\

\section{Le début du spectre des longueurs}

\subsection{Description des premières valeurs}
 Par hypothèse, rappelons que
\[ \frac 12\leq X\leq Y\leq Z<1 \]
Pour ranger dans l'ordre croissant les longueurs des hyperboliques
$\gamma$ vérifiant $\lambda^*(\gamma)\leq 4$, il nous faut
classer les valeurs données dans le tableau ci-contre.\\
\begin{center}
\begin{tabular}{|c|c|c|}
\hline Notation & Formule & Autre écriture\\
\hline$L_0$&$4X^2Y+2XZ-Y$&$\delta(p-1,2)$\\
$L_1$&$2XY+Z$&$\delta(p-1,1)$\\
$L_2$&$2XZ+Y$&$\delta(1,2)$\\
$L_3$&$2YZ+X$&$\delta(2,1)$\\
$L_4$&$4Y^2X+2YZ-X$&$\delta(2,r-1)$\\
$L_5$&$4X^2Z+2XY-Z$&$\delta(1,3)$\\
$L_6$&$8X^3Y+4X^2Z-4XY-Z$&$\delta(1,r-3)$\\
$L_7$&$4Y^2Z+2XY-Z$&$\delta(3,1)$\\
$L_8$&$8Y^3X+4Y^2Z-4XY-Z$&$\delta(3,r-1)$\\
$L_{9}$&$8X^2Y^2+4XYZ-2X^2-2Y^2+1$&$\delta(2,r-2)$\\
$L_{10}$&$4XYZ+2Z^2+2Y^2-1$&$\delta'(2,2)$\\
$L_{11}$&$4XYZ+2X^2+2Y^2+2Z^2-1$&\\
$L_{12}$&$4XYZ+2X^2+2Z^2-1$&\\
$L_{13}$&$4XYZ+2X^2+2Y^2-1$&$\delta(2,2)$\\
\hline
\end{tabular}
\end{center}

Il restera ensuite à examiner les différentes classes de conjugaison
rencontrées dans l'étude précédente pour les valeurs $L_1,L_2,L_3$,
qui se trouveront être les plus petites, de manière à déterminer la
multiplicité de ces longueurs dans le spectre.

Commençons par constater que
\[ 1<L_1\leq L_2\leq L_3 \]
avec $L_1=L_2$ si et seulement si $r=3$ ou $p=q$.

On montre ensuite par des calculs élémentaires que toutes les valeurs, dès qu'elles sont strictement supérieures à l'unité, sont supérieures ou égales à $L_3$, à l'exception des cas suivants:

\begin{align*}&L_0=L_2~~\textrm{ si $r=4$}\\
&L_5=L_1~~\textrm{ si $r=4$}\\
&L_6=L_2~~\textrm{ si $r=5$}\\
&L_7=L_1~~\textrm{ si $p=4$}\\
\end{align*}

Ainsi, on peut énoncer le résultat suivant:\\

\begin{prop}\label{p_13}
On suppose $r\geq 3$. Parmi les éléments hyperboliques de
$\Gamma(r,p,q)$ vérifiant $\lambda^*(\gamma)\leq 4$, les trois plus
petites distances de translation sont
\[ l_1\leq l_2\leq l_3 \]
avec
\begin{align*}
l_1&=2~\mathrm{Argch}~[2\cos\frac{\pi}p\cos\frac{\pi}r+\cos\frac{\pi}q]\\
l_2&=2~\mathrm{Argch}~[2\cos\frac{\pi}q\cos\frac{\pi}r+\cos\frac{\pi}p]\\
l_3&=2~\mathrm{Argch}~[2\cos\frac{\pi}p\cos\frac{\pi}q+\cos\frac{\pi}r]
\end{align*}
\end{prop}

Il reste à  déterminer la multiplicité exacte des longueurs
rencontrées.

\subsection{Etude de la multiplicité}

Dans ce paragraphe, on dit qu'un élément hyperbolique $\gamma$ est
\textbf{associé} à $\delta(k,k')$ si $\gamma$ se décompose en
produit de deux réflexions dont les axes sont disposés comme il est précisé
dans la définition de la quantité $\delta(k,k')$ donnée dans la
première partie de l'article.

Examinons de plus près la multiplicité des longueurs $l_1$ et $l_2$.

Pour cela, nous allons commencer par établir un fait utile dans la
suite:

\begin{lemme}
Si $\gamma$ est associé à $\delta(k,k')$ et $\gamma'$ est associé à
$\delta(p-k,r-k')$, alors $\gamma'$ est conjugué à $\gamma^{-1}$
dans $\Gamma(r,p,q)$.
\end{lemme}

Il suffit pour s'en convaincre de conjuguer $\gamma'$ par un
elliptique d'angle $-2k\pi/p$ bien choisi.\\

\underline{\textit{Multiplicité de $l_1$ si $r\neq p\neq q$ et $r\neq3$}}\\

Reprenons pas à pas l'étude des hyperboliques de niveau $\leq 4$.

Nous constatons que la longueur $l_1$ est atteinte par deux éléments
hyperboliques $\gamma_1,\gamma_2$ associés au cas $(r,p)$ et par
deux éléments associés au cas $(p,r)$. Comme nous l'avons déjà
remarqué en considérant la symétrie des cas, ces derniers sont
conjugués dans $\Gamma(r,p,q)$ aux inverses $\gamma_1^{-1}$ et
$\gamma_2^{-1}$.

On a également trouvé deux éléments $\gamma_3,\gamma_4$ associés au
cas $(r,r,p,r)$ qui possédaient la distance de translation $l_1$ et
deux éléments associés au cas $(r,p,r,r)$, qui sont nécessairement
conjugués aux inverses de $\gamma_3$ et de $\gamma_4$.

On dispose enfin de deux éléments $\gamma_5,\gamma_6$ associés au
cas $(p,p,r,p)$ et de deux éléments associés au cas $(p,r,p,p)$ qui
ont une distance de translation égale à $l_1$, ces derniers étant
conjugués aux inverses de $\gamma_5$ et de $\gamma_6$.\\

Classons maintenant les éléments rencontrés en classe de conjugaison
dans $\Gamma(r,p,q)$. En se servant de décompositions explicites en produits de réflexions bien choisies, on
constate que $\gamma_4$ et $\gamma_5$ sont conjuguées à $\gamma_2$,
que $\gamma_2^{-1}$ est conjugué à $\gamma_3$ puis que $\gamma_3$
est conjugué à $\gamma_6$.

En utilisant le fait établi en préambule, on montre que $\gamma_3$
est conjugué à $\gamma_1$. Par conséquent, on dispose de deux
classes de conjugaison associées à $l_1$: celle de $\gamma_1$ et
celle de $\gamma_1^{-1}$. Comme l'axe de $\gamma_1$ ne passe pas par
un point du pavage, ces éléments ne sont pas conjugués dans
$\Gamma(r,p,q)$ et nous avons bien deux classes de conjugaison
distinctes.\\

Parallèlement à cette étude, on constate que $l_1$ est également
atteinte dans les cas suivants:
\[L_7=L_1 ~~\textrm{ si $p=4$}~~;~~L_5=L_1 ~~\textrm{ si $r=4$}\]
Les éléments qui correspondent à ces valeurs se trouvent cependant à
chaque fois conjugués à $\gamma_1$ ou à $\gamma_1^{-1}$ car
$\delta(3,1)=\delta(p-1,1)$ dans le premier cas et
$\delta(1,3)=\delta(1,r-1)$ dans le second cas.

En conclusion, il n'y a que deux classes de conjugaison dans
$\Gamma(r,p,q)$ qui font intervenir la longueur $l_1$, et celle-ci
est donc de multiplicité 2 dans le spectre des groupes étudiés.\\

\underline{\textit{Multiplicité de $l_2$ si $r\neq p\neq q$ et $r\neq3$}}\\

La longueur $l_2$ est atteinte par deux éléments hyperboliques
$\gamma_1,\gamma_2$ associés au cas $(r,p,r)$, par deux isométries
$\gamma_3,\gamma_4$ associées au cas $(r,r,p)$ et par deux
isométries $\gamma_5$ et $\gamma_6$ correspondant au cas $(p,r,r)$
qui sont nécessairement, par symétrie des types, conjuguées dans
$\Gamma(r,p,q)$ aux inverses
de $\gamma_3$ et de $\gamma_4$.\\

Classons les éléments $\gamma_1,\gamma_2,\gamma_3,\gamma_4$ à
conjugaison près dans $\Gamma(r,p,q)$.

En décomposant judicieusement les éléments étudiés en produits de réflexions, on observe que
$\gamma_3$ est conjugué à $\gamma_2^{-1}$ et que $\gamma_4$ est
conjugué à $\gamma_1^{-1}$. En utilisant maintenant le lemme établi
en début de paragraphe, on montre que $\gamma_2$ est conjugué à
$\gamma_1^{-1}$. Par conséquent, on dispose là encore de deux
classes de conjugaison associées à $l_2$: celle de $\gamma_1$ et
celle de $\gamma_1^{-1}$. Comme l'axe de $\gamma_1$ ne passe pas par
un point du pavage, ces éléments ne sont pas conjugués dans
$\Gamma(r,p,q)$ et nous avons bien deux classes de
conjugaison distinctes.\\

Parallèlement, on observe que $l_2$ est également atteinte dans les
cas suivants:
\[L_6=L_2 ~~\textrm{ si $r=5$}~~;~~L_0=L_2 ~~\textrm{ si $r=4$}\]
Les éléments qui correspondent à ces valeurs se trouvent cependant à
chaque fois conjugués à $\gamma_1$ ou à $\gamma_1^{-1}$ car
$L_6=\delta(1,2)$ dans le premier cas et $L_0=\delta(1,2)$ dans le
second cas.

En conclusion, il n'y a que deux classes de conjugaison dans
$\Gamma(r,p,q)$ qui font intervenir la longueur $l_2$, et celle-ci
est donc de multiplicité 2 dans le spectre des groupes étudiés.\\

\underline{\textit{Multiplicité de $l_2$ et $l_1$ si $p=q$ ou si $r=3$}}\\

Dans ces cas, $l_1=l_2$.\\

Si $r=3$, les éléments associés à ces deux longueurs sont conjugués,
ce qui implique que la multiplicité de la longueur $l_1=l_2$ est
égale à $2$. En effet, $L_1=\delta(1,r-1)=L_2$. Ainsi, si $\gamma_1$
est associé à $l_1$, alors $\gamma_1^{-1}$
est associé à $l_2$ d'après le lemme préliminaire.\\

Traitons maintenant le cas $p=q$. On constate alors, en utilisant des décompositions en produits de refléxions adaptées, que si $\gamma_1$ est associé
à $l_1$ et $\gamma_2$ est associé à
$l_2$ les images des axes de $\gamma_1$ et de $\gamma_2$ sont
différents dans $\mathbf{H}/\Gamma(r,p,q)$ et les deux éléments ne
peuvent donc pas être conjugués.

Ceci montre également que $\gamma_1$ ne peut pas non plus être
conjugué à $\gamma_2^{-1}$.

La multiplicité de la longueur $l_1=l_2$ est alors au moins égale à
4 dans le spectre des groupes considérés. Elle est exactement égale
à 4 si $r\neq p$ car alors la longueur $l_3$ est strictement
plus grande que $l_2$.\\

Si $r=p=q$, cette multiplicité est au moins égale à 5 car, si on
considère $\gamma_3$ un élément associé à $l_3$ correspondant
$\delta(2,1)$, celui-ci n'est pas conjugué aux éléments précédents ni
à leur inverse étant donné que l'axe de $\gamma_3$ est distinct au quotient
des axes de $\gamma_1$ et $\gamma_2$.\\

\underline{\textit{Multiplicité de $l_2$ et $l_1$ si $p=r\neq q$ et $r\neq3$}}\\

Dans ce cas, $l_1<l_2=l_3$ et la multiplicité de $l_1$ est donc
encore
égale à 2.\\

Soit $\gamma_2,\gamma_3$ des éléments associés à $l_2,l_3$ et qui
correspondent à $\delta(1,2)$ et $\delta(2,1)$: d'après ce qui
précède ils sont uniques à conjugaison près et à inverse près.
On constate cependant, là encore en faisant un choix judicieux de représentants, que les
axes sont distincts au quotient: l'élément $\gamma_3$ n'est donc conjugué à aucun des éléments
donnant naissance à la longueur $l_2$. On en déduit donc que la
longueur $l_2=l_3$ a une multiplicité au moins égale à 3 dans les
groupes considérés.

\subsection{Conclusion}\label{s4_4}

On montre ici le résultat principal de ce chapitre

\begin{theoreme}\label{t_2}
 Soit $r\geq 3$. Le spectre des
longueurs du groupe
$\Gamma(r,p,q)$ commence de la manière suivante:\\
$\mathrm{Lsp}~\Gamma(3,3,q)=\{l_1=l_1\dots \}$ \textrm{ pour tout $q\geq 4$}\\
$\mathrm{Lsp}~\Gamma(3,4,4)=\{l_1=l_1\dots \}$\\
$\mathrm{Lsp}~\Gamma(3,p,q)=\{l_1=l_1<l_3\dots \}$ \textrm{ pour tout $p\geq 4,q\geq5$}\\
$\mathrm{Lsp}~\Gamma(r,r,r)=\{l_1=l_1=l_2=l_2=l_3\dots \}$ \textrm{ pour tout $r\geq 4$}\\
$\mathrm{Lsp}~\Gamma(r,r,p)=\{l_1=l_1<l_2=l_2=l_3\dots \}$ \textrm{ pour tout $r\geq4,p\neq r$}\\
$\mathrm{Lsp}~\Gamma(r,p,p)=\{l_1=l_1=l_2=l_2<l_3\dots \}$ \textrm{ pour tout $r\geq4,p\neq r$}\\
$\mathrm{Lsp}~\Gamma(r,p,q)=\{l_1=l_1<l_2=l_2<l_3\dots \}$ \textrm{ pour tout $r\geq4,p\neq r,q\neq p$ }\\
\end{theoreme}

avec les valeurs:
\begin{align*}
l_1&=2~\textrm{Argch}~[2\cos\frac{\pi}p\cos\frac{\pi}r+\cos\frac{\pi}q]\\
l_2&=2~\textrm{Argch}~[2\cos\frac{\pi}q\cos\frac{\pi}r+\cos\frac{\pi}p]\\
l_3&=2~\textrm{Argch}~[2\cos\frac{\pi}p\cos\frac{\pi}q+\cos\frac{\pi}r]
\end{align*}

\textbf{Démonstration:}\\
Supposons dans un premier temps que $\Gamma(r,p,q)$ soit différent
des groupes $\Gamma(3,4,4)$, $\Gamma(4,4,4)$ ou
$\Gamma(5,5,5)$\\

\underline{\textit{Le cas général}}\\
La preuve est alors une utilisation des sections précédentes. Posons
$l_0=l_3$ et décrivons toutes les longueurs inférieures ou égale à
$l_0$.

Si $\gamma$ est un élément hyperbolique de $\Gamma(r,p,q)$, on
examine $\lambda^*(\gamma)$. Si cette valeur est $\geq 5$, alors
tout sommet $x\in\mathcal{E}^*$ vérifie
\[
d(x,\gamma x)\geq\rho^*(D^*(x,\gamma x))\geq
\rho^*(\lambda^*(\gamma))\geq \rho^*(5)
\]
ceci résultant de la croissance de $\rho^*$. On a démontré
précédemment  que $\rho^*(5)>C^*(l_0)$, ce qui impose
\[ \forall
x\in\mathcal{E}^* ~;~d(x,\gamma x)> C^*(l_0) \] La proposition
\ref{p_1} implique que $l(\gamma)>l_0$.

Les éléments hyperboliques de $\Gamma(r,p,q)$ ayant une distance de
translation inférieure ou égale à $l_0=l_3$ possédent donc
nécessairement une distance combinatoire minimale
$\lambda^*(\gamma)\leq 4$. Le théorème résulte alors de la
proposition \ref{p_13} et de l'étude
faite de la multiplicité des longueurs étudiées.\\

Traitons maintenant les autres cas.\\

\underline{\textit{Les cas exceptionnels}}\\
Nous allons déterminer la systole des groupes considérés.

Commençons par remarquer un fait général: si $r\geq 3$ et $2r\geq
p$, $\Gamma(r,p,p)$ est un sous-groupe d'indice $2$ de
$\Gamma(2,p,2r)$.

Ainsi, la systole de $\Gamma(r,p,p)$ est toujours supérieure ou
égale à celle de $\Gamma(2,p,2r)$. On a donc
\[
l_3\geq\textrm{Syst}~\Gamma(r,p,p)\geq\textrm{Syst}~\Gamma(2,p,2r)\]
Rappelons que dans le spectre de $\Gamma(2,p,q)$
(\cite{Philippe2}) nous avons trouvé les longueurs
\begin{align*}
&l_2(1,q-1)=2\textrm{Argch}~[2(\cos\frac{\pi}p)^2+2(\cos\frac{\pi}q)^2-1]\\
&l_2(1,2)=2~\textrm{Argch}~[\cos\frac{\pi}q(4(\cos\frac{\pi}p)^2-1)]\\
&l_1(2)=2~\textrm{Argch}~[2\cos\frac{\pi}p\cos\frac{\pi}q]\\
&l_1(3)=2~\textrm{Argch}~[\cos\frac{\pi}p(4(\cos\frac{\pi}q)^2-1)]
\end{align*} Remarquons enfin que, dans $\Gamma(r,p,p)$,
$l_3=l_2(1,2r-1)$. Ainsi, si $\gamma\in\Gamma(r,p,p)$ vérifie
$l(\gamma)\leq l_3$, alors $\gamma\in\Gamma(2,p,2r)$ et
$l(\gamma)\leq l_2(1,2r-1)$. L'étude faite du spectre des longueurs
des groupes $\Gamma(2,p,q)$ dans \cite{Philippe2} permet d'établir
une liste exhaustive des valeurs possibles pour $l(\gamma)$ en
précisant de plus la classe de conjugaison de l'isométrie $\gamma$
réalisant chaque longueur, ce qui va permettre de conclure.
Reprenons maintenant un à un les cas à étudier.\\

Cherchons la systole du groupe $\Gamma(3,4,4)$.

Si $\gamma\in\Gamma(3,4,4)$ est de longueur $<l_3$, alors
$\gamma\in\Gamma(2,4,6)$ et possède une longueur $< l_2(1,5)$.
D'après l'étude du spectre de $\Gamma(2,4,6)$, il n'existe que deux
solutions: $l(\gamma)=l_1(2)$ ou $l(\gamma)=l_1(3)$. On constate que
$l_1(3)=l_1$, donc cette valeur est atteinte par un élément de
$\Gamma(3,4,4)$. Reste le cas où $l(\gamma)=l_1(2)$. On constate
alors que $\gamma$ ne laisse pas stable $\mathcal{E}^*$ donc
$\gamma$ ne peut pas appartenir à $\Gamma(3,4,4)$: la seule
valeur du spectre de $\Gamma(3,4,4)$ qui est $<l_3$ est donc $l_1=l_2$.\\
Examinons la question du groupe $\Gamma(4,4,4)$. Cette fois-ci,
$\Gamma(4,4,4)$ est contenu dans $\Gamma(2,4,8)$ et la seule valeur
de ce spectre $<l_3=l_2(1,7)$ est $l_1(2)$. Celle-ci correspond là
encore à un élément de $\Gamma(2,4,8)$ qui ne préserve pas
$\mathcal{E}^*$ et qui n'appartient donc pas à $\Gamma(4,4,4)$.

Reste enfin le cas du groupe $\Gamma(5,5,5)$. Soit $\gamma$ un
élément de $\Gamma(5,5,5)$ avec $l(\gamma)<l_3$. L'élément $\gamma$
est alors contenu dans $\Gamma(2,5,10)$, et il n'existe que deux
valeurs $<l_3=l_2(1,9)$ dans le spectre des longueurs de ce groupe.
Ces valeurs sont $l_1(2)$ et $l_2(1,2)$ (on rappelle que dans le
groupe considéré ici, $l_2(1,q-1)=l_1(3)$).

Dans ces deux cas, $\gamma$ ne préserve pas $\mathcal{E}^*$ et ne
peut donc pas appartenir à $\Gamma(5,5,5)$. Par conséquent, dans les
cas $\Gamma(4,4,4)$ et $\Gamma(5,5,5)$, la systole est $l_3$ et sa
multiplicité est donnée par l'étude déjà effectuée au cas général:
elle est $\geq5$. Dans le cas $\Gamma(3,4,4)$, la systole est $l_1$
et sa multiplicité est $\geq2$.\\
Ceci achève la preuve du théorème. $\maltese$

\section{Compléments sur le spectre des longueurs de certains groupes}
Pour des raisons qui apparaîtront naturellement lorsque l'on voudra différencier deux groupes de triangles à l'aide de leur spectre, nous allons dans cette section affiner l'étude du spectre pour certains groupes de triangles que la systole seule ne permettra pas de caractériser.

\subsection{Compléments et retour sur le spectre avec $r=2$}

Reprenons les notations de \cite{Philippe2}, où l'on travil avec $\mathcal{P}$ le sous pavage de $\mathcal{P}_0$ obtenu
en considérant uniquement les sommets de valence $q$ et les arêtes
les reliant: celui-ci est constitué de $p$-gones de côté $2c$ avec

\[\cosh c =\frac{\cos\frac{\pi}p}{\sin\frac{\pi}q}
\]

et les angles aux sommets sont tous égaux à $2\pi/q$.

Si $x,y$ sont deux sommets de $\mathcal{P}$, on considère l'ensemble
des chemins $\beta$ du pavage formés d'arêtes consécutives de
$\mathcal{P}$ et reliant $x$ et $y$. Le nombre d'arêtes d'un tel
chemin est sa longueur, notée $L(\beta)$. La distance combinatoire
entre $x$ et $y$ est
\[ D(x,y)=\textrm{inf}\{ L(\beta)~;~\beta\textrm{ reliant $x$ et $y$
}  \} \] Fixons maintenant $s_0$ dans $\mathcal{P}$. Pour $n\geq 1$,
on pose
\[ \rho(n)=\mathrm{inf}\{d(s_0,s)~;~D(s_0,s)=n\} \]
Si $l_0>0$, on introduit la constante
\[ C(l_0)=\mathrm{argch}~[(\cosh c)^2(\cosh l_0 -1)+1] \]
Enfin, pour tout $\gamma$ élément hyperbolique de $\Gamma(2,p,q)$,
on définit
\[ \lambda(\gamma)=\textrm{inf}\{~D(s,\gamma s)~;~s\in\mathcal{P}\} \]
qui mesure le déplacement minimal de $\gamma$ sur $\mathcal{P}$ au
sens combinatoire du terme.

Nous avons montré dans \cite{Philippe2} que $\rho$ est une fonction
croissante et que, pour chaque groupe $\Gamma(2,p,q)$ fixé, si l'on cherche à décrire les valeurs du spectre inférieure à une valeur $l_0$ il
suffit de trouver $n_0$ avec $\rho(n_0)>C(l_0)$ puis de décrire les
éléments hyperboliques de niveau $< n_0$. Dans la mise en pratique,
nous avons déterminé les longueurs des éléments hyperboliques de
niveau $\leq 2$ pour chacun
des groupes, ce qui a permi d'établir les résultats suivants:\\

$\mathrm{Lsp}~\Gamma(2,3,7)=\{l_2(1,q-1)<\dots \}$\\
$\mathrm{Lsp}~\Gamma(2,3,q)=\{l_2(1,q-1)<l_1(4)\dots \}$ \textrm{ pour tout $q\geq 8$}\\
$\mathrm{Lsp}~\Gamma(2,4,q)=\{l_1(2)=l_1(2)<l_1(3)=\dots<l_2(1,q-1)<\dots\}$ \textrm{ pour $q=6,7$}\\
$\mathrm{Lsp}~\Gamma(2,4,q)=\{l_1(2)=l_1(2)<l_2(1,q-1)<\dots \}$ \textrm{ pour tout $q=5,q\geq8$}\\
$\mathrm{Lsp}~\Gamma(2,5,5)=\{l_1(2)=l_1(2)<l_2(1,q-1)<\dots \}$ \\
$\mathrm{Lsp}~\Gamma(2,5,q)=\{l_1(2)=l_1(2)<l_2(1,2)\dots \}$ pour tout $q\geq 6$\\
$\mathrm{Lsp}~\Gamma(2,p,q)=\{l_1(2)=l_1(2)<\dots \}$ \textrm{pour tout $p\in\{6,7,8,9,10\}$ }\\
$\mathrm{Lsp}~\Gamma(2,p,q)=\{l_1(2)=l_1(2)<l_2(1,2)\dots \}$ \textrm{ pour tout $p\geq 11$}\\

avec
\begin{align*}
&l_2(1,q-1)=2~\textrm{Argch}~[2(\cos\frac{\pi}p)^2+2(\cos\frac{\pi}q)^2-1]\\
&l_2(1,2)=2~\textrm{Argch}~[\cos\frac{\pi}q(4(\cos\frac{\pi}p)^2-1)]\\
&l_1(2)=2~\textrm{Argch}~[2\cos\frac{\pi}p\cos\frac{\pi}q]\\
&l_1(3)=2~\textrm{Argch}~[\cos\frac{\pi}p(4(\cos\frac{\pi}q)^2-1)]\\
\end{align*}
Remarquons que les résultats énoncés ici sont plus précis que ceux donnés
dans \cite{Philippe2}: nous avons en effet déterminé la
multiplicité exacte de la systole. Expliquons comment.\\

\subsubsection{Détermination des multiplicités}

Nous reprenons ici les notations introduites dans \cite{Philippe2}
pour décrire les éléments de déplacement combinatoire $\leq 2$, à savoir :
\begin{align*}
l_2(k,k')&=2\mathrm{Argch}~\mid
\sin\frac{k'\pi}q\sin\frac{k\pi}q\cosh 2c
-\cos\frac{k\pi}q\cos\frac{k'\pi}q \mid\\
l_1(k)&= 2~\mathrm{Argch}~[\sin\frac{k\pi}q
\frac{\cos\frac{\pi}p}{\sin\frac{\pi}q}]
\end{align*}

Commençons par énoncer le

\begin{lemme}
On a $l_1(k)=l_1(q-k)$ et les éléments associés sont, à conjugaison près
dans $\Gamma(2,p,q)$, inverses l'un de l'autre.

De même, $l_2(k,k')=l_2(q-k',q-k)$ et les éléments associés sont, là encore à
conjugaison près dans $\Gamma(2,p,q)$, inverses l'un de l'autre.
\end{lemme}

Pour s'en convaincre, on utilise la conjugaison par un retournement
en ayant choisi des représentants adéquats pour
chacune de ces longueurs.\\

\underline{\textit{La multiplicité de $l_1(2)$.}}\\

D'après l'étude des hyperboliques de niveau $\leq2$ effectuée dans
\cite{Philippe2}, la longueur $l_1(2)$ ne peut être égalée que dans
les cas suivants
\begin{displaymath}
\left\{\begin{array}{l}
l_1(2)=l_1(q-2)~~;~~\forall q\geq4\\
l_1(2)=l_1(3)~~\textrm{ si $q=5$ }\\
l_1(2)=l_2(1,2)~~\textrm{ si $p=q=5$ }
\end{array}\right.
\end{displaymath}

Remarquons également que pour n'importe quelle valeur de $q$, on a
\[ l_2(1,2)=l_2(2,1)=l_2(q-2,q-1)=l_2(q-1,q-2) \]

En introduisant des représentants bien choisis des éléments réalisant la longueur $l_1(2)$ et en s'aidant du résultat énoncé précédemment, on constate que ces éléments sont soit conjugués, soit inverses l'un de l'autre.

Comme l'axe des hyperboliques considérés ne passe pas par un point du
pavage initial la multiplicité de $l_1(2)$ est exactement égale à 2 dans tous les groupes considérés.\\

\underline{\textit{La multiplicité de $l_2(1,q-1)$.}}\\

On a vu dans l'étude des hyperboliques de niveau $\leq 2$ effectuée
dans \cite{Philippe2} que $l_2(1,q-1)$ ne peut être égalée que dans
les situations suivantes:

\begin{displaymath}
\left\{\begin{array}{l}
l_2(1,q-1)=l_2(q-1,1)~~;~~\forall~q\geq 4\\
l_2(1,q-1)=l_1(3)~~\textrm{ si $p=3$}
\end{array}\right.
\end{displaymath}

De manière analogue, on constate que ces éléments sont soit conjugués, soit inverses l'un de l'autre.
Cependant, l'axe de chacun passe par un point de valence paire
du pavage et la longueur
$l_2(1,q-1)$ a une multiplicité exactement égale à 1 dans tous les groupes
étudiés.\\

\subsubsection{La seconde longueur de certains groupes
$\Gamma(2,p,q)$}

Nous allons maintenant décrire la deuxième longueur des groupes
$(2,p,\infty)$ si $6\leq p \leq 10$; ainsi que la
deuxième longueur de $\Gamma(2,10,10)$.\\

Si $p\geq6$, nous avons montré dans \cite{Philippe2} que
\[
\rho(3)=\frac{16Y^2(Y^2+Z^2-1)(2Y^2-1)+2Y^2+Z^2-1}{1-Z^2}\] Nous
avons également établi que
\[C(l_2(1,q-1))=\frac{8Y^2(Y^2+Z^2-1)(Z^2+Y^2)+1-Z^2}{1-Z^2}\]
Un nouveau calcul montre que
\[C(l_2(1,2)) =\frac{2Y^2(4ZY^2-Z-1)(4ZY^2-Z+1)+1-Z^2}{1-Z^2}
\]
Si $6\leq p\leq 10$ et si $q=\infty$, on constate que
\[\rho(3)-C(l_2(1,2))=16Y^4(2Y^2-1)+2Y^2 -16Y^4(2Y^2-1)=2Y^2 >0 \]
Ceci montre que
\[ \fbox{$\textrm{Lsp}~\Gamma(2,p,\infty)=\{ l_1(2)=l_1(2)<l_2(1,2)\dots \}~~;~~6\leq p\leq 10$} \]
Sachant que, pour tous ces groupes, on a constaté dans
\cite{Philippe2} que la seule longueur inférieure à $l_2(1,2)$ parmi
les
éléments de niveau $1,2$ est $l_1(2)$.\\

Etudions maintenant le cas du groupe $\Gamma(2,10,10)$. On constate
que, si $p=q=10$,
\[\rho(3)-C(l_2(1,q-1))>0\]

On en déduit donc comme précédemment que
\[\fbox{$\mathrm{Lsp}~\Gamma(2,10,10)=\{l_1(2)=l_1(2)<l_2(1,2)\dots
\}$}\]

\subsection{Compléments sur le spectre avec $r=3$}

Nous allons déterminer ici la deuxième longueur du spectre des groupes $\Gamma(3,3,5),~\Gamma(3,3,6)$ et $\Gamma(3,4,4)$.\\

\subsubsection{Calcul de $\rho^*(5)$ pour les groupes $\Gamma(3,3,q)$}
Examinons dans la démarche effectuée dans les sections précédentes ce qu'il advient
si $r=p=3$. Les seuls types donnant réellement un chemin de longueur
combinatoire $5$ sont $(r,p,r,p,r)$ et $(p,r,p,r,p)$. Dans le cas $(r,p,r,p,r)$, on a constaté que le sommet $y$ le plus proche vérifie
\[ \cosh
x_0y=\frac{(X^2+Y^2+Z^2+2XYZ-1)[2Z(1+Y+Z+XZ)+(2Z^2-1)(1+X)]^2}{(1+X)(1+Y)(1+Z)}+1
\]
Par conséquent, dans les groupes $\Gamma(3,3,q)$:
\[\cosh\rho^*(5)\geq16Z^5+8Z^4-12Z^3-2Z^2+3Z+\frac12 \]
avec égalité si $q\geq6$.\\

\subsubsection{Choix d'un nouveau $l_0$ pour les groupes $\Gamma(3,3,q)$}

Nous aurons besoin de différencier les spectres de
$\Gamma(3,3,5)$ et $\Gamma(2,5,5)$ d'un côté et de $\Gamma(3,3,6)$
et $\Gamma(2,4,12)$ d'un autre côté. Nous constaterons que ces
groupes possèdent la même systole, et il nous faudra considérer
les deuxièmes longueurs. Choisissons donc pour $l_0$ la deuxième valeur
du spectre des groupes $\Gamma(2,5,5)$ (si $q=5$) et
$\Gamma(2,4,12)$ (si $q=6$).

\begin{center}
\begin{tabular}{|c|c|c|}
\hline Valeur de $q$ & Choix de $l_0$ & Valeur de $C^*(l_0)$\\
\hline $5$&$2\textrm{Argch}~[4(\cos\frac{\pi}5)^2-1]$&$\displaystyle{\frac{365\sqrt{5}+833}{324}}$\\
$6$&$2\textrm{Argch}~[2(\cos\frac{\pi}{12})^2+2(\cos\frac{\pi}4)^2-1]$&$\displaystyle{\frac{785}{324}\sqrt{3}+\frac{61}{18}}$\\
\hline
\end{tabular}
\end{center}

Pour les deux valeurs de $q$ étudiées, la comparaison des valeurs de
$\cosh\rho^*(5)$ et de $\cosh C^*(l_0)$ permet d'établir le fait
suivant:

\begin{lemme}
La deuxième longueur du spectre des groupes $\Gamma(3,3,q)~;~q=5,6$
est égale à
\[
l_2=2~\mathrm{Argch}~[2(\cos\frac{\pi}q)^2+\cos\frac{\pi}q-\frac12]
\]
et cette longueur a une multiplicité supérieure ou égale à $2$.
\end{lemme}

\textbf{Démonstration:}\\
On vérifie que dans les deux cas
$\cosh\rho^*(5)>\cosh C^*(l_0)$ ce qui permet d'établir que tous les
hyperboliques de $\Gamma(3,3,5)$ et $\Gamma(3,3,6)$ de niveau
supérieur à $5$ ont une longueur strictement supérieure à $l_0$. On
examine maintenant les longueurs obtenues précédemment pour
les hyperboliques de niveau inférieur à $4$ (cf. le tableau suivant)
pour constater que la deuxième longueur de $\Gamma(3,5,5)$ et
$\Gamma(3,3,6)$ est aussi égale à $l_0$, obtenue avec
$L_{10}=L_{12}$. La multiplicité est $\geq2$ car l'axe ne passe pas
par un point de valence paire.$\maltese$\\

\begin{center}
\begin{tabular}{|c|c|}
\hline Notation & Valeur dans le groupe $\Gamma(3,3,q)$\\
\hline$L_0,L_4,L_{13}$&$Z$\\
$L_1,L_2,L_3,L_9$&$Z+\frac12$\\
$L_5,L_7$&$\frac12$\\
$L_6,L_8$&$-\frac12$\\
$L_{10},L_{12}$&$2Z^2+Z-\frac12$\\
$L_{11}$&$2Z^2+Z$\\
\hline
\end{tabular}
\end{center}
Dans ce tableau, nous ne retenons que les valeurs $>1$ et la
longueur $l$ associée est déterminée par $\cosh(l/2)=L$.\\

Examinons maintenant le cas du groupe $\Gamma(3,4,4)$ qui est, rappelons le, un sous-groupe d'indice 2 dans
$\Gamma(2,4,6)$. Or nous connaissons le début du spectre de
$\Gamma(2,4,6)$:
\[\mathrm{Lsp}~\Gamma(2,4,6)=\{l_1(2)=l_1(2)<l_1(3)=\dots<l_2(1,5)\dots\}
\]
Nous avons constaté dans l'étude précédente que $l_1=l_1(3),
l_3=l_2(1,5)$ et que ces deux valeurs étaient donc réalisées par des
éléments de $\Gamma(3,4,4)$, contrairement à $l_1(2)$. Par
conséquent, il n'existe aucune valeur entre $l_1$ et $l_3$ dans le
spectre de $\Gamma(3,4,4)$:
\[\fbox{$\mathrm{Lsp}~\Gamma(3,4,4)=\{l_1=l_1=\dots<l_3\dots \}$}\]


\begin{thebibliography}{99}

\bibitem{Bavard} C.~Bavard. \emph{L'aire systolique conforme des groupes christallographiques du plan},~Ann. Inst. Fourier,~\textbf{43}(1993),~\no 3,~p.815-842.
\bibitem{Beardon} A.F.~Beardon.\emph{The geometry of Discrete Groups},~Graduate Texts in Mathematics
    ,~\textbf{91},~Springer-Verlag,~1983.
\bibitem{Buser2} P.~Buser,~K.-D.~Semmler. \emph{The geometry and spectrum of the one holed
torus},~Comment. Math. Helv.,~\textbf{63} (1988),~p.259-274.
\bibitem{Dal'bo} F.~Dal'Bo. \emph{Trajectoires géodésiques et
horocycliques},~Savoirs Actuels,~CNRS Editions,~2007.
\bibitem{Dianu} R.~Dianu. \emph{Sur le spectre des tores pointés},
Thèse,~ EPFL, Lausanne, 2000.
\bibitem{Haas} A.~Haas. \emph{Length spectra as moduli for hyperbolics surfaces},~Duke Math. J.,~\textbf{52} (1985),~p.922-935.
\bibitem{Kac}M.~Kac. \emph{Can One Hear the Shape of a Drum?},~The American Mathematical Monthly,~\textbf{73},~\no 4, Part 2: Papers
in Analysis (Apr., 1966),~p. 1-23.
\bibitem{LW}R.~Lehman,~C.~White. \emph{Hyperbolic billiards
path}.\\\textsf{http://www.tilings.org/pubs/TRlehmanwhite.pdf}
\bibitem{Philippe1} E.~Philippe. \emph{Géométrie des surfaces hyperboliques},
Thèse,~ Université Paul Sabatier, Toulouse, 2008.
\bibitem{Philippe2} E.~Philippe. \emph{Les groupes de triangles $(2,p,q)$ sont déterminés par leur spectre des longueurs},
~Ann. Inst. Fourier,~\textbf{58}(2008),~\no 7,~p.2659-2693.\\\textsf{http://arxiv.org/abs/0807.4746}
\bibitem{Ratcliffe} J.G.~Ratcliffe. \emph{Foundations of Hyperbolic Manifolds},~Graduate Texts in Mathematics
    149,~Springer-Verlag,~1994.
\bibitem{Vogeler} R.~Vogeler. \emph{On the geometry of Hurwitz surfaces},~Thesis,~Univ. Florida, 2003.
\end{thebibliography}
\end{document}